\documentclass[12pt]{amsart}

\usepackage{dsfont}
\usepackage{amsmath}
\usepackage{amssymb}
\usepackage{graphicx}

\usepackage{amssymb}

\usepackage{amsfonts}

\usepackage{soul}

\usepackage{mathpazo}
\usepackage{color}
\usepackage{yfonts}

\usepackage{paralist}

\usepackage{stmaryrd}

\usepackage{amsxtra}

\def\carre{\hfill $\Box$}

\def\th{\theta}

\def\cA{\mathcal A}
\def\eps{\varepsilon}

\newcommand{\norm}[1]{\bigl\| #1 \bigr\|}
\renewcommand{\norm}[1]{|\!|\!| #1 |\!|\!|}

\def\a{          \alpha}

\def\cA{          \mathcal A}

\def\th{        \theta}

\def \R{{\mathbb R}}

\def \Z{{\mathbb Z}}
\def \N{{\mathbb N}}

\newcommand{\T}{{\mathbb T}}
\newcommand{\prf}{{\begin{proof}}}
\newcommand{\epf}{{\end{proof}}}

\newcommand{\Q}{{\mathbb Q}}

\theoremstyle{definition}

\def\bee{\begin{equation}}
\def\eee{\end{equation}}

\theoremstyle{remark}

\newcommand{\pdvr}[2]
{\dfrac{\partial^{#2} #1}{\partial \theta^{#2_1} \partial r^{#2_2}}}
\newcommand{\pdvrs}[2]
{\partial^{#2} #1 /\partial \theta^{#2_1} \partial r^{#2_2}}
\newtheorem{thm}{Theorem}
\newtheorem{corr}{Corollary}

\newtheorem{alemma}{Lemma}
\newtheorem{definition}{Definition}
\numberwithin{equation}{section}

\begin{document}

\title{On the convergence to $0$ of $m_n�\xi \ {mod \ }   1$  }
\author[Bassam Fayad and Jean-Paul Thouvenot]{Bassam Fayad and Jean-Paul Thouvenot}

\address{
(Fayad) IMJ-PRG CNRS \\ (Thouvenot) LPMA UPMC CNRS}
\email{bassam@math.jussieu.fr, jean-paul.thouvenot@upmc.fr}

\date{\today}

\begin{abstract} We show that for any irrational number $\a$ and a sequence of integers $\{m_l\}_{l\in \N}$ such that  $\displaystyle{\lim_{l\to \infty} \norm{m_l \a} = 0}$, there exists a continuous measure $\mu$ on the circle such that $\displaystyle{\lim_{l\to \infty}  \int_\T \norm{m_l \th} d\mu(\th) = 0}$. This implies that any rigidity sequence of any ergodic transformation is a rigidity sequence for some weakly mixing dynamical system. 

On the other hand, we show that for any $\a \in \R - \Q$, there exists a sequence of integers $\{m_l\}_{l\in \N}$ such that  $\norm{m_l \a} \to 0$  and $m_l \theta [1]$ is dense on the circle if and only if $\th \notin \Q \a+\Q$.
\end{abstract}

\maketitle

\section{Introduction}
We denote $\T=\R/\Z$. For $x \in \R$ we denote by $\norm{x}:=\inf_{n\in \N} |x-n|$. We denote by $x [1]$ the fractional part of $x$. 

\medskip 

In this note, we prove the following results. 

\begin{thm} \label{th2} For any $\a \in \R - \Q$ and a sequence of integers $\{m_l\}_{l\in \N}$ such that  $\lim_{l\to \infty} \norm{m_l \a} = 0$, there exists a measure $\mu$ on $\T$ which has no atoms such that $\lim_{l\to \infty}  \int_\T \norm{m_l \th} d\mu(\th) = 0$.
\end{thm}

  \begin{thm} \label{th1} For any $\a \in \R - \Q$, there exists a sequence of integers $\{m_l\}_{l\in \N}$ such that  $\norm{m_l \a} \to 0$  and  such that $m_l \theta [1]$ is dense in $\T$ if and only if $\th \notin \Q \a+\Q$.
\end{thm}

Due to the Gaussian measure space construction (see \cite{CFS}, or for example Proposition 2.30 in \cite{lem}),  Theorem \ref{th2} has a direct consequence on rigidity sequences of weakly mixing dynamical systems.  
The following is in fact an equivalent statement to Theorem \ref{th2}. 
\begin{corr} \label{cor2} For any $\a \in \R- \Q$ and a sequence of integers $\{m_l\}_{l\in \N}$ such that  $\lim_{l\to \infty} \norm{m_l \a} = 0$, there exists a weak mixing dynamical system $(T,M,m)$ such that $\{m_l\}_{l\in \N}$ is a rigidity sequence for $(T,M,m)$. 
\end{corr}

A consequence of Corollary \ref{cor2} is a positive answer to a question raised in \cite{lem}, namely :  {\it a rigidity sequence of any ergodic transformation (on a probability space without atoms) with discrete spectrum is a rigidity sequence for some weakly mixing dynamical system}. Indeed, Corollary \ref{cor2}  deals with the case of a pure point spectrum with an irrational rotation of the circle as a factor. The case of a purely rational spectrum was treated in \cite{lem},  Proposition 3.27. In case the spectrum is purely rational, our proof of Theorem \ref{th2} given below applies with only one modification: instead of working with the orbit of $0$ under the rotation $R_{\alpha}$, ($\alpha \notin Q/Z)$, one considers the union of the orbits of $0$ under the actions of all the finite groups which appear in the (necessarily dense) support of the spectral measure. 

 A completely different  solution to the same question was given in \cite{adams}, who proved directly Corollary \ref{cor2} based on a sophisticated and involved cut and stack construction.


In contrast, our proof is much simpler and is based on the straightforward characterization of rigidity as a spectral property, which reduces the answer to the question to the construction of a continuous probability measure on the circle with a Fourier transform converging to $1$ along the rigidity subsequence as stated in Theorem \ref{th2}. This possible approach to the question was discussed in detail in \cite{lem}.

The second result, Theorem \ref{th1}, asserts that it is not possible to expect more than what is obtained in Theorem \ref{th2}, namely that a strong convergence of $\norm{m_n \th}$ to $1$ on an uncountable set $K$ is not possible in general for a sequence $m_n$ such that $\lim_{l\to \infty} \norm{m_l \a} = 0$, $\a \in \R-\Q$. Constructing such a set $K$ was a possible strategy to proving Corollary \ref{cor2} (see for example Proposition 3.3 in \cite{lem}), and Theorem \ref{th1} shows that this approach cannot be adopted in general. 

Given an increasing sequence of integers $m_n$, the study of the accumulation points on the circle of the sequence $\{m_n \xi \}$, for $\xi$ irrational, has a long history and a  rich literature (see for example \cite{bugeaud,dub} and references therein). 
Weyl \cite{weyl} proved that for any increasing sequence $m_n$, it holds that for almost every $\xi$, $\{ m_n \xi \}$ is dense on the circle. The set of irrationals $\xi$ such that $\{ m_n \xi \}$ is not dense in $\T$ is called the set of exceptional points for the sequence $m_n$. Our result asserts the existence for any $\a \in \R-\Q$ of a sequence $m_n$ for which the set of exceptional points is reduced to $\Q \a+\Q$. To our knowledge, no other examples of increasing sequences $m_n$ with a countable exceptional set are known in the literature.

\medskip {\sc Acknoweledgment.} The authors are grateful to the referee for suggesting improvements to the first version of the paper.

  \section{Proof of Theorem \ref{th2}}

Fix $\a \in \R-\Q$  and a sequence of integers $\{m_l\}_{l\in \N}$ such that  
$$\hspace{1.3cm} (\star) \ \ \ \hspace{1cm}  \lim_{l\to \infty} \norm{m_l \a} = 0.\hspace{2cm} $$

For a probability measure $\mu$ on $\T$ we write $\mu^n= |\int_\T \norm{m_n \th} d\mu(\th)|$.

We will construct a sequence $\mu_p$, $p\geq 0$,  of probability measures on $\T$ of the form $\frac{1}{2^p} \sum_{i=1}^{2^p} \delta_{x_i}$ with $x_i=k_i \a$ such that there exists an increasing  sequence $\{N_p\}$ for which
\begin{itemize}
\item[(1)] For every $p\geq 1$, for every $j \in [0,p-1]$, for every $n \in [N_{j},N_{j+1}]$, $\mu_p^n<\frac{1}{2^{j}}$
\item[(2)] For  every $p_0 \in \N^*$, if we let 
$$\eta_{p_0}=\frac{1}{4} \inf_{1\leq i<i'\leq 2^{p_0}}   \norm{k_i\a-k_{i'}\a}$$
 then for every $ l  \in \N$ and every $r \in [1,2^{p_0}]$, $\norm{k_{l2^{p_0} +r} \a-k_r\a} < \eta_{p_0}$
 \item[(3)] $\mu_p^n < \frac{1}{2^{p+1}}$ for $n\geq N_p$. 
\end{itemize} 
When going from the measure $\mu_p$ to $\mu_{p+1}$ we will {\it add} $2^p$ masses at points selected nearby $x_1,\ldots,x_{2^p}$ that are already chosen for $\mu_p$. 

Theorem \ref{th2} clearly follows from the above construction. Indeed, property (1) will imply that any weak limit $\mu_\infty$ of $\mu_p$ satisfies $\mu_\infty^n \to 0$. While by (2) we get that for each $p_0$ the intervals $(k_r \a-\eta_{p_0},k_r \a+\eta_{p_0})$, $r \in [1,2^{p_0}]$, on the circle are disjoint and each have mass $1/2^{p_0}$ for every $\mu_{p}$, $p\geq p_0$, and hence for $\mu_\infty$ that has therefore no atoms. 

Property (3) is not necessary in the proof of the theorem, but it is useful to fulfill the inductive hypotheses (1) and (2) of the construction.

For $p=0$, we let $k_1=0$ and $\mu_0$ is thus the Dirac measure at $0$. We let $N_0=0$. 
For $p=1$, we let $k_2=1$ so $\mu_1$ is the average of the Dirac measures at $0$ and at $\alpha$. Observe that for any $n$, $\mu_1^n<\frac{1}{2}$ which fulfills (1) for $p=1$. 
We also choose $N_1$ sufficiently large so that  $\mu_1^n < \frac{1}{2^{2}}$ for $n\geq N_1$, thelatter being possible due to $(\star)$.

We now assume selected $k_i$ for $i\leq 2^p$  and $N_l$ for $l\leq p$ such that (1) and (3) are satisfied up to $p$, and (2) is satisfied for every $p_0\leq p$ and every $0\leq l\leq 2^{p-p_0}-1$.

We choose $k_{2^p+1}$ such that $k_{2^p+1}\a$ is sufficiently close to $k_1 \a$ so that 
$$\nu_{p,1}=\frac{1}{2^{p+1}} \sum_{i=1}^{2^{p+1}} \delta_{k'_i\a}$$
where $k'_i=k_i$ for $i\leq 2^p$ and $k'_{2^p+1}=k_{2^p+1}$ while $k'_{2^p+r}= k_{r}$ for $r\in [2,2^p]$
satisfies $\nu_{p,1}^n < \frac{1}{2^{j}}$ for every $n \in [N_{j},N_{j+1}]$ and $j\in  [0, p-1]$. 

Since for every $n$ we have that $|\nu_{p,1}^n-\mu_p^n|< \frac{1}{2^{p+1}} \norm{m_n k_{2^p+1} \a - m_n k_{1} \a} < \frac{1}{2^{p+1}}$ we get by (3) that $\nu_{p,1}^n< \frac{1}{2^{p+1}}+ \frac{1}{2^{p+1}} < \frac{1}{2^{p}}$ for every $n\geq N_p$. 
Next we choose, $N_{p,1}>N_p$ sufficiently large such that $\nu_{p,1}^n< \frac{1}{2^{p+2}} $ for $n \geq N_{p,1}$, which is possible by $(\star)$. In this way, we  select inductively $k_{2^p+s}$, then  $N_{p,s}$ for $s=1,2,\ldots, 2^p$, and  set
$$\nu_{p,s} = \frac{1}{2^{p+1}} \sum_{i=1}^{2^{p+1}} \delta_{k'_i\a}$$ 
where $k'_i=k_i$ for $i\leq 2^p+s$ and $k'_{2^p+t}=k_{t}$  for $t\in [s+1,2^p]$. Choosing for each $s$, $k_{2^p+s} \a$ sufficiently close to $k_s \a$ then $N_{p,s}$ sufficiently large we can insure that 
\begin{itemize}
\item $\nu_{p,s}^n < \frac{1}{2^{j}}$ for every $n \in [N_{j},N_{j+1}]$ and $j\leq p-1$. 
\item $\nu_{p,s}^n < \frac{1}{2^{p}}$ for every $n \geq N_{p}$
\item $\nu_{p,s}^n < \frac{1}{2^{p+2}}$ for every $n \geq N_{p,s}$
\end{itemize}
The first point can be established inductively due to the fact that if  $k_{2^p+s}\a$ is chosen very close to $k_s \a$   then the measures $\nu_{p,s-1}$ and $\nu_{p,s}$ are very close. 

The same argument gives the second point for $N_p\leq n\leq N_{p,s-1}$. As for $n \geq N_{p,s-1}$ we use the fact that $|\nu_{p,s}^n-\nu_{p,s-1}^n|< \frac{1}{2^{p+1}}$ and the fact that $\nu_{p,s-1}^n < \frac{1}{2^{p+2}}$ for every $n \geq N_{p,s-1}$ to conclude that $\nu_{p,s}^n< \frac{1}{2^{p}}$. For the third point we just choose $N_{p,s}$ sufficiently large and use $(\star)$.

Finally, we let $N_{p+1}=N_{p,2^p}$ and $\mu_{p+1}=\nu_{p,2^p}$ and observe that the measure $\mu_{p+1}$ satisfies (1).

Also, since $k_{2^p+s} \a$ can be chosen arbitrarily close to  $k_s \a$ for $s=1,\ldots, 2^p$, we get that
for every $p_0\leq p+1$, for every $l= 2^{p-p_0}+l'-1$, $l'\leq 2^{p-p_0}$, $\norm{k_{l2^{p_0} +r} \a} \sim \norm{k_{l'2^{p_0} +r} \a} \sim 
 \norm{k_{r} \a}$ from where (2) follows for $p+1$. The proof of Theorem \ref{th2} is thus complete. \carre

   \section{Proof of Theorem \ref{th1}}

In all this section $\a \in \R-\Q$ is fixed.

\begin{definition} For an interval $I \subset \T$, for $\eps>0$, and integers $N_1<N_2$, we say that $\th \in \cA(N_1,N_2,I,\eps,\a)$ if for  every $m \in [N_1,N_2)$ such that $\norm{m\a} <\eps$ we have that $\{ m \th \}\notin I$. 
\end{definition}

\begin{alemma} \label{lemma1} For every $l\geq 2$, there exists $L(l)\in \N$, such that for every $0<\eps\leq\frac{1}{2l^2}$, for every  $\nu>0$, $N \in \N$, there exist $K(\eps)>0$ and $N'=N'(l,\eps,\nu,N)\in \N$ such  that $\th \in \cA(N,N',I,\eps,\a)$ for some interval $I$ of size $1/l$, implies that $\norm{k\a-s\th}<\nu$ for some $|k|\leq K(\eps)$ and some $|s|\leq L(l)$. 
\end{alemma}
\medskip 

\noindent {\it Proof.} For any $\eps>0$, consider an approximation by trigonometric polynomials of  $2 \chi_{\eps}$ where $\chi_{\eps}$ is the characteristic function of the subset of $\T$ $ [0,\eps]\cup [1-\eps,1]$, namely $\phi_\eps : \T \to \R$ such that 
\begin{itemize}
\item $\phi_\eps(x) > 1$ for every $x \in [0,\eps]\cup [1-\eps,1]$
\item $\phi_\eps(x)>-\eps^3$ for every $x \in \T$ 
\item There exists $K\in \N$ such that $\phi_\eps(x) = \sum_{|k|\leq K} \hat{\phi}_k e^{i2 \pi kx}$.
\end{itemize} 

Similarly, for $l \geq 2$, let $\varphi_l: \T \to \R$ be such that 
\begin{itemize}
\item $\varphi_l(y) > 1$ for every $y \notin [0,\frac{1}{l}]$
\item $|\varphi_l(y)| < l^2$ for every $y \in \T$
\item There exists $L\in \N$ such that $\varphi_l(y) = \sum_{0<|k|\leq L} \hat{\varphi}_k e^{i2 \pi ky}$.
\end{itemize} 
Note that the second requirement includes the fact that $\int \varphi_l(y) dy =0$. 

For $\psi : \T^2 \to \R$ and $(\a,\th) \in \R^2$ we define for $k \in \N$
$$S^{\a,\th}_{k}\psi(x,y)= \sum_{i=0}^{k-1} \psi(x+i\a,y+i\th).$$

Fix $I=[y_0,y_0+\frac{1}{l}]$, for some $y_0 \in \T, l \geq 2$.

Define $\psi_{\eps,l} : \T^2 \to \R$ by $\psi_{\eps,l}(x,y)= \phi_\eps(x) \varphi_l(y-y_0)$. For $N' \in \N$ sufficiently large there exist more than $\eps^2N'$ integers  $i\in [N,N')$ such that $\norm{i\a} <\eps$. If $\th \in \cA(N,N',I,\eps,\a)$  then $S^{\a,\th}_{N'}  \psi_{\eps,l} (0,0) > (\eps^2 - l^2 \eps^3)N' \geq \frac{1}{2} \eps^2 N'$. 

On the other hand, we have that  
\begin{equation*} S^{\a,\th}_{N'} \psi_{\eps,l} (x,y)= \sum_{|k|\leq K, 0<|j|<L} \hat{\phi}_k \hat{\varphi}_j \frac{1-e^{i2 \pi N' (k\a+j\th)}}{1-e^{i2 \pi (k\a+j\th)}} e^{i2 \pi (kx+jy)} \end{equation*}
hence, if $\norm{k\a-j\th} \geq \nu$ for every $|k|\leq K$ and every $0<|j|\leq L$,  then  $S^{\a,\th}_{N'} \psi_{\eps,l}(x,y)$ is bounded  independently of $N'$ which contradicts $S^{\a,\th}_{N'} \psi_{\eps,l}(0,0) > \frac{1}{2} \eps^2 N'$.  

\carre 

\medskip

\noindent {\it Proof of Theorem \ref{th1}.} 
 Define for $n\geq 1$,  the sequence $l_n=n+1$ and $L_n:=L(l_n)$ given by Lemma \ref{lemma1}. Let  $\eps_n=\frac{1}{2 (n+1)^2}$  and $K_n=K(\eps_n)$ given by Lemma \ref{lemma1}. Let $\nu_n=\frac{1}{n} \inf_{0<|k|\leq (n+1) K_{n+1}} \norm{k\a}$. Take $N_0=0$ and apply Lemma \ref{lemma1} with $l=l_1$, $\eps=\eps_1$, $N=N_0$ and $\nu=\nu_1$. Denote then $N_1=N'(l_1,\eps_1,\nu_1,N_0)$. We then apply inductively Lemma \ref{lemma1} with  $l=l_n$, $\eps=\eps_n$, $N=N_n$ and $\nu=\nu_n$ and choose $N_{n+1}$ arbitrarily large such that $N_{n+1}\geq N'(l_n,\eps_n,\nu_n,N_n)$. 

We define an increasing sequence $m_l$ by taking successively  for every $i$, all the integers  $m \in [N_{i},N_{i+1})$ such that 
$\norm{m\a} < \eps_i$ (choosing $N_{n+1}$ to be sufficiently large in our inductive construction guarantees that the sequence $m_n$ is not empty). 

Suppose now $\th$ is such that $m_n \th [1] $ is not dense on the circle. Then, there exists $k$ and an interval $I$ of size $l_k$ such that for every $n$, $m_n \th [1]  \notin I$. In other words, $\th \in \cA(N_n,N_{n+1},I,\eps_n,\a)$ for every $n\geq n_0$, $n_0$ sufficiently large.  Let $L=L_k$.   By Lemma \ref{lemma1} we get that  $\norm{k_n\a-l\th} <\nu_n$  for some $|k_n|\leq K_n$ and some $0<|l|\leq L$. Hence $\norm{k'_n \a-L! \th} <L! \nu_n$  for some $|k'_n|\leq L! K_n$. It follows that $\norm{(k'_{n+1}-k'_n)\a}<2 L! \nu_n$. From the definition of $\nu_n$ this implies that $k'_{n+1}=k'_n$ for sufficiently large $n$, say $n\geq n_1$. Since $\nu_n\to 0$, we get that $\norm{k'_{n_1} \a -L! \th }=0$, which gives $\th \in \Q \a +\Q$. Conversely, $\{ m_n \th \}$ for $\th \in \Q \a +\Q$ is clearly not dense on the circle and Theorem \ref{th1} is proved. \carre

\end{document}